\documentclass[a4paper]{amsart}
\newtheorem{theorem}{Theorem}
\newtheorem{definition}{Definition}
\newtheorem{proposition}{Proposition}
\def\End{\mathop{\mathrm{End}}}
\def\Ad{\mathop{\mathrm{Ad}}}
\def\id{id}
\def\vsp{V}

\def\la{b}
\def\cla{c_b}

\def\SURFACE{\Sigma}
\def\NTR{n_\Sigma}
\def\VERTICES{V}

\def\TRIANGLES{T}
\def\REALS{\mathbb R}
\def\COMPLEXS{\mathbb C}
\def\INTEGERS{\mathbb Z}
\def\PGROUP{\mathbb S}
\def\SDIT{\Delta_\SURFACE}

\def\a{\mathsf a}
\def\x{\mathsf x}
\def\u{\mathsf b}
\def\e{\mathsf e}
\def\E{\mathsf E}
\def\f{\mathsf f}
\def\g{\mathsf g}
\def\RMAT{\mathsf R}
\def\BRAID{\mathsf B}
\def\LENGTH{\mathsf L}
\def\ROTATE{\mathsf A}
\def\PTOLEMY{\mathsf T}
\def\PERMUTE{\mathsf P}
\def\MOM{\mathsf p}
\def\POS{\mathsf q}
\def\QDILOG{\mathop{\mathrm{e}_b}}
\def\IMUN{\mathsf i}
\def\FUNCTOR{\mathsf F}
\def\DEHN{\mathsf D}
\def\CONSTRAINT{\mathsf z}
\def\MCG{\mathcal M_\SURFACE}

\title[On the spectrum of Dehn twists]{On 
the spectrum of Dehn twists in quantum Teichm\"uller theory}
\author{R.~M.~Kashaev}
\address{Steklov Mathematical Institute at St. Petersburg,
Fontanka 27, St. Petersburg 191011, Russia}
\curraddr{Helsinki Institute of Physics, P.O. Box 9, (Siltavuorenpenger
  20C) FIN-00014
University of Helsinki, Finland}
\email{kashaev@pdmi.ras.ru}
\date{August 2000}
\begin{document}
\begin{abstract}
  The operator realizing a Dehn twist in quantum Teichm\"uller theory
  is diagonalized and continuous spectrum is obtained. This result is
  in agreement with the expected spectrum of conformal weights in
  quantum Liouville theory at $c>1$. The completeness condition of the
  eigenvectors includes the integration measure which appeared in the
  representation theoretic approach to quantum Liouville theory by
  Ponsot and Teschner. The underlying quantum group structure is also
  revealed.
\end{abstract}
\maketitle
\section{Introduction}
\label{sec:introduction}

The quantization problem of Teichm\"uller spaces of punctured surfaces
is connected to quantum Chern--Simons theory with non-compact gauge
groups (and thereby to 2+1-dimensional quantum gravity \cite{ach,wit})
and non-rational quantum conformal field theory in two dimensions. In
particular, quantum Teichm\"uller theory is expected to describe the
space of conformal blocks in quantum Liouville theory \cite{ver}. From
the mathematical point of view this can be useful in three-dimensional
topology, in particular for construction of new link and
three-manifold invariants.
    
The Teichm\"uller spaces of punctured surfaces have been quantized in
two different but essentially equivalent ways: as (degenerate) Poisson
manifolds \cite{fock2,fock1}, and as simplectic manifolds with
Weil--Petersson simplectic structure \cite{kash1}. The both approaches
employ the real analytic coordinates closely related to the Penner
parameterization of the decorated Teichm\"uller spaces \cite{penn}.
The main outcome of the theory is the construction of a particular
projective (infinite dimensional) representation of the mapping class
groups of punctured surfaces. The projective factor has been
shown to be related to the Virasoro central charge in quantum
Liouville theory \cite{kash2}.

In this paper we describe the spectrum of Dehn twists in the quantum
Teichm\"uller theory, derive an explicit formula for the associated
with braidings $R$-matrix (square of which is a Dehn twist), and
uncover the underlying infinite dimensional representations of the
quantum group $\mathcal U_q(\mathfrak{sl}(2))$. Those representations
were considered in papers \cite{schm,pon1,pon2}.
 
The paper is organized as follows.  In
section~\ref{sec:quantum-dilogarithm} we briefly remind the main
properties of the non-compact quantum dilogarithm (closely related to
the double sine function, see for example \cite{jim}). This is the
main technical object entering practically all the constructions of
the paper.  In section~\ref{sec:basic-algebr-syst} we describe the
algebraic system underlying the quantum Teichm\"uller theory and its
realization in terms of the quantum dilogarithm.
Section~\ref{sec:mapping-class-group} is a survey of the part of
papers \cite{kash1,kash2} which describes the mapping class group
representation in terms of the basic algebraic system of
section~\ref{sec:basic-algebr-syst}. Solution of the Dehn twist
spectral problem is described in section~\ref{sec:diag-dehn-twist}
(proof of Theorem~\ref{theor:spec} is to appear in a separate
publication).  Section~\ref{sec:braiding-r-matrix} contains derivation
of the $R$-matrix and revealing the associated quantum group
structure.

\emph{Acknowledgments.} I wish to thank L.~D.~Faddeev for valuable
discussions and encouragement. This work is supported by RFBR grant
99-01-00101, INTAS grant 99-01705, and Finnish Academy.

\section{Quantum dilogarithm}
\label{sec:quantum-dilogarithm}

Let complex $\la$ have a nonzero real part $\Re \la\ne0$.  The
\emph{non-compact quantum dilogarithm}, $\QDILOG(z)$, is defined by
the integral formula~\cite{fad95}
   \begin{equation}\label{eq:ncqdl}
\QDILOG(z)\equiv\exp\left(\frac{1}{4}
\int_{\IMUN 0-\infty}^{\IMUN 0+\infty}
\frac{e^{-\IMUN 2 zw}\, dw}{\sinh(w\la)
\sinh(w/\la) w}\right)
    \end{equation}
    in the strip $|\Im z|<|\Im\cla|$, where
\[
\cla\equiv\IMUN(\la+\la^{-1})/2.
\]
In particular, when $\Im\la^2>0$, we have the following product
formula:
\begin{equation}\label{eq:ratio}
\QDILOG(z)=
(e^{2\pi (z+\cla)\la};q^2)_\infty/
(e^{2\pi (z-\cla)\la^{-1}};\bar q^2)_\infty,
\end{equation}
where
\[
q\equiv e^{\IMUN\pi\la^2},\quad\bar q\equiv e^{-\IMUN\pi\la^{-2}}.
\]
Using the symmetry properties
\[
\QDILOG(z)=\mathop{\mathrm{e}_{-b}}(z)=\mathop{\mathrm{e}_{1/b}}(z),
\]
we assume in what follows
\[
\Re\la>0,\quad \Im \la\ge0.
\]
Function $\QDILOG(z)$ can be analytically continued in variable $z$ to
the entire complex plane as a meromorphic function with essential
singularity at infinity, and with the following characteristic
properties:
\begin{description}
\item[poles and zeros]
  \begin{equation}
    \label{eq:polzer}
(\QDILOG(z))^{\pm1}=0\ \Leftrightarrow \ z=
\mp(\cla+m\IMUN\la+n\IMUN\la^{-1}),\ m,n\in\INTEGERS_{\ge0};
  \end{equation}
\item[behavior at infinity] depending on the direction along which the
  limit is taken,
\begin{equation}\label{eq:asymp}
\QDILOG(z)\bigg\vert_{|z|\to\infty} \approx\left\{
\begin{array}{ll}
1&|\arg z|>\frac{\pi}{2}+\arg \la;\\
e^{\IMUN\pi z^2-\IMUN\pi (1+2\cla^2)/6}&|\arg z|<\frac{\pi}{2}-\arg\la;\\
\frac{(\bar q^2;\bar q^2)_\infty}{\Theta(\IMUN\la^{-1}z;-\la^{-2})}&
|\arg z-\frac\pi2|<\arg\la;\\
\frac{\Theta(\IMUN\la z;\la^{2})}{(q^2; q^2)_\infty}&
|\arg z+\frac\pi2|<\arg\la,
\end{array}\right.
\end{equation}
where
\[
\Theta(z;\tau)\equiv\sum_{n\in\INTEGERS}
e^{\IMUN\pi\tau n^2+\IMUN2\pi zn}, \quad\Im\tau>0;
\]
\item[inversion relation]
\begin{equation}\label{eq:inversion}
\QDILOG(z)\QDILOG(-z)=e^{\IMUN\pi z^2
-\IMUN\pi(1+2\cla^2)/6};
\end{equation} 
\item[functional equations]
\begin{equation}\label{eq:shift}
\QDILOG(z-\IMUN\la^{\pm1}/2)=(1+e^{2\pi\la^{\pm1}z })
\QDILOG(z+\IMUN\la^{\pm1}/2);
\end{equation}
\item[unitarity] when $\la$ is either real or on unit circle,
\begin{equation}\label{eq:qdlunitarity}
(1-|\la|)\Im\la=0\ \Rightarrow\ 
\overline{\QDILOG(z)}=1/\QDILOG(\bar z) 
;
\end{equation}
\item[pentagon equation]
\begin{equation}\label{eq:pent}
\QDILOG(\MOM)\QDILOG(\POS)=\QDILOG(\POS)\QDILOG(\MOM+\POS)
\QDILOG(\MOM),
\end{equation}
if self-adjoint operators $\MOM$ and $\POS$ in $L^2(\REALS)$ satisfy
the Heisenberg commutation relation $[\MOM,\POS]=(2\pi \IMUN)^{-1}$;
\item[integral analogue of Ramanujan's summation formula]
\begin{multline}\label{eq:raman}
  \int_{\REALS}
  \frac{\QDILOG(x+u)}{\QDILOG(x+v)}e^{2\pi\IMUN wx}\, dx\\
  = \frac{\QDILOG(u-v-\cla)\QDILOG(w+\cla)}{\QDILOG(u-v+w-\cla)}
  e^{-2\pi\IMUN w(v+\cla)+\IMUN\pi(1-4\cla^2)/12}
  \\
  =\frac{\QDILOG(v-u-w+\cla)}{\QDILOG(v-u+\cla)\QDILOG(-w-\cla)}
  e^{-2\pi\IMUN w(u-\cla)-\IMUN\pi(1-4\cla^2)/12},
\end{multline}
where
\begin{equation}\label{eq:restrictions1}
\Im(v+\cla)>0,\quad\Im(-u+\cla)>0, \quad \Im(v-u)<\Im w<0.
\end{equation}
\end{description}
Restrictions (\ref{eq:restrictions1}) can actually be relaxed by
deforming the integration path in the complex $x$ plane, keeping the
asymptotic directions of the two ends within the sectors $\pm(|\arg
x|-\pi/2)>\arg\la$. The enlarged in this way domain for the variables
$u,v,w$ in eqn~\eqref{eq:raman} has the form:
\begin{equation}\label{eq:restrictions2}
|\arg (\IMUN z)|<\pi-\arg\la,\quad z\in\{w,v-u-w,u-v-2\cla\}.
\end{equation}
As the matter of fact the pentagon equation~\eqref{eq:pent} is
equivalent to the integral Ramanujan formula~\eqref{eq:raman}, see
 \cite{fkv} for the proof.

\section{The basic algebraic system}
\label{sec:basic-algebr-syst}

Introduce an important notation.
Let $\vsp$ be a vector space.  For any natural $1\le i\le m$ we define
embeddings
        \[
        \iota_i\colon \End\vsp\ni \a\mapsto
        \a_i=\underbrace{\id\otimes\cdots\otimes\id}_{i-1\ 
          \mathrm{times}}\otimes \a\otimes\id\otimes\cdots\otimes \id
        \in\End\vsp^{\otimes m} ,
         \]
         ie $\a$ stands in the $i$-th position.
         If $\u\in \End\vsp^{\otimes k}$ for some $1\le k\le m$ and
         $\{i_1,i_2,\ldots,i_k\}\subset\{1,2,\ldots,m\}$, we write
    \[
\u_{i_1i_2\ldots i_2}\equiv\iota_{i_1}\otimes\iota_{i_2}\otimes\cdots
\otimes\iota_{i_k}(\u).
    \]
Note also that the permutation group $\PGROUP_m$
is naturally represented in $\vsp^{\otimes m}$:
    \begin{equation}\label{eq:perm}
\PERMUTE_\sigma (x_1\otimes\cdots\otimes
x_i\otimes\cdots) =x_{\sigma^{-1}(1)}\otimes\cdots\otimes
x_{\sigma^{-1}(i)}\otimes\ldots,\quad\sigma\in \PGROUP_m.
    \end{equation}
    
The projective representation of the mapping class groups of punctured
surfaces, associated with the quantum Teichm\"uller theory, is based
on the concrete realization of the following algebraic system of
equations on two invertible elements $\ROTATE\in \End\vsp$ and
$\PTOLEMY\in\End\vsp^{\otimes 2}$,

\begin{gather}
  \label{eq:1}
\ROTATE^3=\id, \\
 \label{eq:2}
\PTOLEMY_{12}\PTOLEMY_{13}\PTOLEMY_{23}=\PTOLEMY_{23}\PTOLEMY_{12}, \\
 \label{eq:3}
\ROTATE_1\PTOLEMY_{12}\ROTATE_2=\ROTATE_2\PTOLEMY_{21}\ROTATE_1, \\
  \label{eq:4}
\PTOLEMY_{12} \ROTATE_1\PTOLEMY_{21}=\zeta\ROTATE_1\ROTATE_2\PERMUTE_{(12)}, 
\end{gather}
where $\zeta$ is a (nonzero) complex number, $\PERMUTE_{(12)}$ being
defined in eqn~\eqref{eq:perm} with 
\[
\PGROUP_2\ni\sigma=(12)\colon 1\mapsto2\mapsto1.
\]

\subsection{Useful notation}
\label{sec:useful-notation}

In practical calculations it is very convenient to use the following
 notation which
permits to avoid writing explicitly operator $\ROTATE$.
For any $\a\in\End V$ we denote
\begin{equation}
  \label{eq:14}
  \a_{\hat k}\equiv \ROTATE_k\a_k\ROTATE_k^{-1},\quad
 \a_{\check k}\equiv \ROTATE_k^{-1}\a_k\ROTATE_k.
\end{equation}
Obviously
\[
\a_{\check{\hat k}}=\a_{\hat{\check k}}=\a_k,\quad \a_{\hat{\hat
    k}}=\a_{\check k},\quad \a_{\check{\check k}}=\a_{\hat k},
\]
where the last two equations follow from eqn~\eqref{eq:1}. Besides, it
is useful to denote
\begin{equation}
  \label{eq:16}
 \PERMUTE_{(kl\ldots m\hat k)}\equiv\ROTATE_k\PERMUTE_{(kl\ldots
  m)},\quad
\PERMUTE_{(kl\ldots m\check k)}\equiv\ROTATE_k^{-1}\PERMUTE_{(kl\ldots
  m)}, 
\end{equation}
where $(kl\ldots m)$ is cyclic permutation,
\[
(kl\ldots m)\colon k\mapsto l\mapsto\ldots\mapsto m\mapsto k.
\]
Eqns~\eqref{eq:3},~\eqref{eq:4} in this notation take very compact
form
\begin{gather}
  \label{eq:15}
  \PTOLEMY_{12}=\PTOLEMY_{\hat2\check1},\\
  \PTOLEMY_{12}\PTOLEMY_{2\hat1}= \zeta\PERMUTE_{(12\hat1)}.
\end{gather}
\subsection{Realization through the quantum dilogarithm}
\label{sec:real-thro-quant}

Take
\[
\vsp=L^2(\REALS).
\]
Let self-adjoint operators $\MOM,\POS$ satisfy the Heisenberg
commutation relation
    \[
    [\MOM,\POS]=(2\pi\IMUN)^{-1}.
    \]
    Then operators
\begin{gather}
  \label{eq:5}
  \ROTATE\equiv e^{-\IMUN\pi/3}e^{\IMUN
    3\pi\POS^2}e^{\IMUN\pi(\MOM+\POS)^2}\in L^2(\REALS),\\
\label{eq:t-in-t-of-psi}
\PTOLEMY\equiv e^{\IMUN 2\pi\MOM_1\POS_2}
(\QDILOG(\POS_1+\MOM_2-\POS_2))^{-1}\in L^2(\REALS^2),
    \end{gather}
    do satisfy equations~\eqref{eq:1}~--~\eqref{eq:4} with
    $\zeta=e^{\IMUN\pi\cla^2/3}$.
    Note that operator $\ROTATE$ here is unitary and 
    is characterized by the equations (up to a normalization factor)
\[
\ROTATE\POS\ROTATE^{-1}=\MOM-\POS,\quad\ROTATE\MOM\ROTATE^{-1}=-\POS,
\]
while operator $\PTOLEMY$ is unitary if $(1-|\la|)\Im\la=0$.

\section{The mapping class group representation}
\label{sec:mapping-class-group}

Here we recall how  the mapping class groups of
punctured surfaces are represented in terms
of system~\eqref{eq:1}~--~\eqref{eq:4}.

\subsection{Decorated ideal triangulations}
\label{sec:decor-ideal-triang}

We call a two-cell in a cell complex (or CW-complex, see \cite{dub})
{\em triangle} if exactly three boundary points of the corresponding
two-disk are mapped to the zero-skeleton.  We shall also call
zero-cells and one-cells {\em vertices} and {\em edges}, respectively.

Let $\SURFACE$ be a compact oriented (possibly with boundary) surface
with finite non-empty set $\VERTICES_\Sigma$ of marked points called
also \emph{punctures}.
\begin{definition}
  A cell complex decomposition of $\SURFACE$ is called \emph{ideal
    triangulation} if
\begin{itemize}
  \item the zero-skeleton coincides with $\VERTICES_\Sigma$;

  \item all two-cells are triangles.
\end{itemize}
Two ideal triangulations are considered \emph{equivalent} if there
exists an isotopy of $\SURFACE$ fixed on
$\partial\SURFACE\cup\VERTICES_\SURFACE$ deforming one into another.
\end{definition}
We suppose that $\SURFACE$ admits ideal triangulations (for example,
each boundary component must contain at least one marked point). Any
ideal triangulation of $\SURFACE$ has one and the same number of
triangles $\NTR$.  Denote $\TRIANGLES(\tau)$ the set of triangles in
ideal triangulation $\tau$.
\begin{definition}
  Ideal triangulation $\tau$ is called {\em decorated} if
  \begin{itemize}
  \item each triangle is provided by a marked corner;
    
  \item all triangles are numbered, ie a bijective {\em numbering}
    mapping
    \[
    \bar\tau\colon \{1,\ldots,\NTR\}\rightarrow\TRIANGLES(\tau)
    \]
    is fixed.
  \end{itemize}
\end{definition}
Graphically (see below) in a triangle we put the corresponding integer
inside of it, and asterisk at the marked corner. The set of all
decorated ideal triangulations of $\SURFACE$ will be denoted $\SDIT$.
For terminological convenience and if no confusion is possible, in
what follows we sometimes will use the term triangulation as a
substitute for decorated ideal triangulation.

\subsubsection{Action of the permutation group}

The permutation group $\PGROUP_{\NTR}$ naturally acts in $\SDIT$ from
the right,
    \[
    \SDIT\times\PGROUP_{\NTR}\ni(\tau,\sigma)\mapsto\tau\sigma
    \in\SDIT,
    \]
    by changing the numbering mapping:
    \[
    \overline{\tau\sigma}=\bar\tau\circ\sigma.
    \]

\subsubsection{Changing of a marked corner}

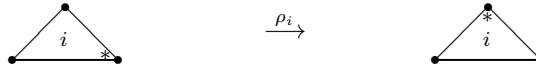
\begin{figure}[htb]
  \centering
\begin{picture}(200,20)
\put(0,0){\begin{picture}(40,20)
\put(0,0){\line(1,0){40}}
\put(0,0){\line(1,1){20}}
\put(20,20){\line(1,-1){20}}
\put(0,0){\circle*{3}}
\put(20,20){\circle*{3}}
\put(40,0){\circle*{3}}
\footnotesize
\put(33,0){$*$}
\put(18,5){$i$}
\end{picture}}
\put(160,0){\begin{picture}(40,20)
\put(0,0){\line(1,0){40}}
\put(0,0){\line(1,1){20}}
\put(20,20){\line(1,-1){20}}
\put(0,0){\circle*{3}}
\put(20,20){\circle*{3}}
\put(40,0){\circle*{3}}
\footnotesize
\put(17.5,14){$*$}
\put(18,5){$i$}
\end{picture}}
\put(95,8){$\stackrel{\rho_i}{\longrightarrow}$}
\end{picture}
\caption{Transformation $\rho_i$ changes the
  marked corner of triangle $\bar\tau(i)$.}\label{ft}
\end{figure}

If $\tau\in\SDIT$ and $1\le
    i\le\NTR$, then triangulation
    $\rho_i(\tau)$ is obtained from $\tau$ by changing the
    marked corner of triangle $\bar\tau(i)$ as is shown in
    figure~\ref{ft}.

\subsubsection{The flip transformation}

\begin{figure}[htb]
  \centering
\begin{picture}(200,40)
\put(0,0){
\begin{picture}(40,40)
\put(20,0){\line(-1,1){20}}
\put(40,20){\line(-1,-1){20}}
\put(0,20){\line(1,1){20}}
\put(40,20){\line(-1,1){20}}
\put(20,0){\line(0,1){40}}
\put(20,0){\circle*{3}}
\put(0,20){\circle*{3}}
\put(20,40){\circle*{3}}
\put(40,20){\circle*{3}}
\footnotesize
\put(10,18){$i$}\put(26,18){$j$}
\put(1,18){$*$}
\put(19.5,2){$*$}
\end{picture}}
\put(160,0){\begin{picture}(40,40)
\put(20,0){\line(-1,1){20}}
\put(40,20){\line(-1,-1){20}}
\put(0,20){\line(1,1){20}}
\put(40,20){\line(-1,1){20}}
\put(0,20){\line(1,0){40}}
\put(20,0){\circle*{3}}
\put(0,20){\circle*{3}}
\put(20,40){\circle*{3}}
\put(40,20){\circle*{3}}
\footnotesize
\put(18,26){$i$}\put(18,10){$j$}
\put(3,20){$*$}
\put(17.5,1){$*$}
\end{picture}}
\put(95,17){$\stackrel{\omega_{ij}}{\longrightarrow}$}
\end{picture}
\caption{The flip transformation $\omega_{ij}$  clockwise ``rotates''
  one diagonal of the quadrilateral until it matches another
  diagonal.}\label{fe}
\end{figure}
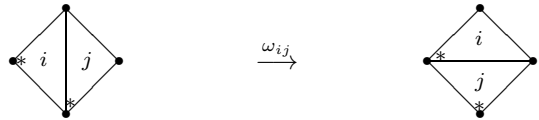

Let two distinct triangles $\bar\tau(i)$, $\bar\tau(j)$ have a common
edge and their marked corners be as in the lhs of figure~\ref{fe},
so the common edge is a diagonal of a quadrilateral combined of the two
triangles.
Then triangulation $\omega_{ij}(\tau)$ is obtained from $\tau$ by replacing
the common edge by the opposite diagonal of the quadrilateral, assigning the numbers and
marked corners to new triangles as is shown in the rhs of
figure~\ref{fe}. Note that this \emph{flip} transformation
$\omega_{ij}$ implicitly depends on the triangulation it transforms as at
fixed $i$ and $j$ it is not defined on all triangulations. 

\subsubsection{The properties}
\label{sec:properties}

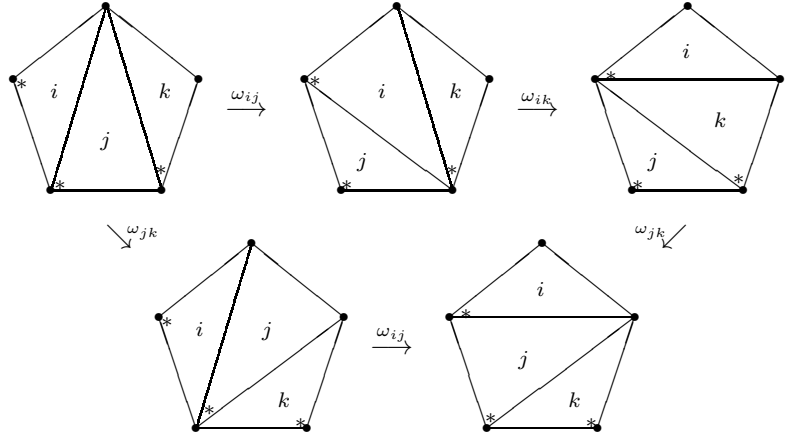
\begin{figure}[htb]
  \centering
\begin{picture}(290,160)
\put(0,90){\begin{picture}(70,70)
\put(14,0){\line(-1,3){14}}
\put(56,0){\line(1,3){14}}
\put(0,42){\line(5,4){35}}
\put(70,42){\line(-5,4){35}}
\put(14,0){\line(1,0){42}}
\qbezier(14,0)(20,20)(35,70)
\qbezier(56,0)(50,20)(35,70)
\put(14,0){\circle*{3}}
\put(56,0){\circle*{3}}
\put(0,42){\circle*{3}}
\put(70,42){\circle*{3}}
\put(35,70){\circle*{3}}
\footnotesize
\put(1,38.5){$*$}
\put(15.5,0){$*$}
\put(53.75,6){$*$}
\put(14,35){$i$}
\put(33,17){$j$}
\put(55,35){$k$}
\end{picture}}
\put(110,90){\begin{picture}(70,70)
\put(14,0){\line(-1,3){14}}
\put(56,0){\line(1,3){14}}
\put(0,42){\line(5,4){35}}
\put(70,42){\line(-5,4){35}}
\put(14,0){\line(1,0){42}}
\put(56,0){\line(-4,3){56}}
\qbezier(56,0)(50,20)(35,70)
\put(14,0){\circle*{3}}
\put(56,0){\circle*{3}}
\put(0,42){\circle*{3}}
\put(70,42){\circle*{3}}
\put(35,70){\circle*{3}}
\footnotesize
\put(2,39.5){$*$}
\put(14,0){$*$}
\put(53.75,6){$*$}
\put(28,35){$i$}
\put(20,9){$j$}
\put(55,35){$k$}
\end{picture}}
\put(220,90){\begin{picture}(70,70)
\put(14,0){\line(-1,3){14}}
\put(56,0){\line(1,3){14}}
\put(0,42){\line(5,4){35}}
\put(70,42){\line(-5,4){35}}
\put(14,0){\line(1,0){42}}
\put(56,0){\line(-4,3){56}}
\put(0,42){\line(1,0){70}}
\put(14,0){\circle*{3}}
\put(56,0){\circle*{3}}
\put(0,42){\circle*{3}}
\put(70,42){\circle*{3}}
\put(35,70){\circle*{3}}
\footnotesize
\put(4,41.5){$*$}
\put(14,0){$*$}
\put(52.25,2.5){$*$}
\put(33,50){$i$}
\put(20,9){$j$}
\put(45,24 ){$k$}
\end{picture}}
\put(55,0){\begin{picture}(70,70)
\put(14,0){\line(-1,3){14}}
\put(56,0){\line(1,3){14}}
\put(0,42){\line(5,4){35}}
\put(70,42){\line(-5,4){35}}
\put(14,0){\line(1,0){42}}
\qbezier(14,0)(20,20)(35,70)
\put(14,0){\line(4,3){56}}
\put(14,0){\circle*{3}}
\put(56,0){\circle*{3}}
\put(0,42){\circle*{3}}
\put(70,42){\circle*{3}}
\put(35,70){\circle*{3}}
\footnotesize
\put(1,38.5){$*$}
\put(17,5){$*$}
\put(51.5,0){$*$}
\put(14,35){$i$}
\put(39,35){$j$}
\put(45,8){$k$}
\end{picture}}
\put(165,0){\begin{picture}(70,70)
\put(14,0){\line(-1,3){14}}
\put(56,0){\line(1,3){14}}
\put(0,42){\line(5,4){35}}
\put(70,42){\line(-5,4){35}}
\put(14,0){\line(1,0){42}}
\put(0,42){\line(1,0){70}}
\put(14,0){\line(4,3){56}}
\put(14,0){\circle*{3}}
\put(56,0){\circle*{3}}
\put(0,42){\circle*{3}}
\put(70,42){\circle*{3}}
\put(35,70){\circle*{3}}
\footnotesize
\put(4,41.5){$*$}
\put(13.5,2){$*$}
\put(51.5,0){$*$}
\put(33,50){$i$}
\put(26,24){$j$}
\put(45,8){$k$}
\end{picture}}
\put(35,70){$\searrow$}\put(43,74){\tiny$\omega_{jk}$}
\put(245,70){$\swarrow$}\put(235,74){\tiny$\omega_{jk}$}
\put(135,28){$\stackrel{\omega_{ij}}{\longrightarrow}$}
\put(80,118){$\stackrel{\omega_{ij}}{\longrightarrow}$}
\put(190,118){$\stackrel{\omega_{ik}}{\longrightarrow}$}
\end{picture}
  \caption{Proof of the pentagon equation~\eqref{eq:22}.}
  \label{fig:pen-om}
\end{figure}

We have the
following properties:
\begin{gather}
  \label{eq:19}
  \rho_i\circ\rho_i\circ\rho_i=id,\\\label{eq:22}
  \omega_{jk}\circ\omega_{ik}\circ\omega_{ij}=\omega_{ij}\circ\omega_{jk},
\\\label{eq:23}
  (\rho_i^{-1}\times\rho_j)\circ\omega_{ij}=\omega_{ji}\circ(\rho_i^{-1}\times\rho_j),\\
\label{eq:20}
\omega_{ji}\circ\rho_i\circ\omega_{ij}=(ij)\circ(\rho_i\times\rho_j).
\end{gather}
The first equation is evident since there are only three possibilities to
mark a corner in a triangle. The other three equations are proved
pictorially in figures~\ref{fig:pen-om} --~\ref{fig:inv-rel}.
\begin{figure}[htb]
  \centering
  \begin{picture}(200,100)
\put(0,60){\begin{picture}(200,40)
\put(0,0){\begin{picture}(40,40)
\put(20,0){\line(-1,1){20}}
\put(40,20){\line(-1,-1){20}}
\put(0,20){\line(1,1){20}}
\put(40,20){\line(-1,1){20}}
\put(20,0){\line(0,1){40}}
\put(20,0){\circle*{3}}
\put(0,20){\circle*{3}}
\put(20,40){\circle*{3}}
\put(40,20){\circle*{3}}
\footnotesize
\put(10,18){$i$}\put(26,18){$j$}
\put(1,18){$*$}
\put(19.5,2){$*$}
\end{picture}}
\put(160,0){\begin{picture}(40,40)
\put(20,0){\line(-1,1){20}}
\put(40,20){\line(-1,-1){20}}
\put(0,20){\line(1,1){20}}
\put(40,20){\line(-1,1){20}}
\put(0,20){\line(1,0){40}}
\put(20,0){\circle*{3}}
\put(0,20){\circle*{3}}
\put(20,40){\circle*{3}}
\put(40,20){\circle*{3}}
\footnotesize
\put(18,26){$i$}\put(18,10){$j$}
\put(3,20){$*$}
\put(17.5,1){$*$}
\end{picture}}
\put(95,17){$\stackrel{\omega_{ij}}{\longrightarrow}$}
\end{picture}}
\put(18,46){$\downarrow$\tiny$\rho_i^{-1}\times\rho_j$}
\put(145,46){{\tiny$\rho_i^{-1}\times\rho_j$}$\downarrow$}
\put(0,0){\begin{picture}(200,40)
\put(0,0){\begin{picture}(40,40)
\put(20,0){\line(-1,1){20}}
\put(40,20){\line(-1,-1){20}}
\put(0,20){\line(1,1){20}}
\put(40,20){\line(-1,1){20}}
\put(20,0){\line(0,1){40}}
\put(20,0){\circle*{3}}
\put(0,20){\circle*{3}}
\put(20,40){\circle*{3}}
\put(40,20){\circle*{3}}
\footnotesize
\put(10,18){$i$}\put(26,18){$j$}
\put(16,32){$*$}
\put(34,17.5){$*$}
\end{picture}}
\put(160,0){\begin{picture}(40,40)
\put(20,0){\line(-1,1){20}}
\put(40,20){\line(-1,-1){20}}
\put(0,20){\line(1,1){20}}
\put(40,20){\line(-1,1){20}}
\put(0,20){\line(1,0){40}}
\put(20,0){\circle*{3}}
\put(0,20){\circle*{3}}
\put(20,40){\circle*{3}}
\put(40,20){\circle*{3}}
\footnotesize
\put(18,26){$i$}\put(18,10){$j$}
\put(18,34){$*$}
\put(32,15){$*$}
\end{picture}}
\put(95,17){$\stackrel{\omega_{ji}}{\longrightarrow}$}
\end{picture}}
\end{picture}
  \caption{Proof of the symmetry relation~\eqref{eq:23}.}
  \label{fig:sym-om}
\end{figure}
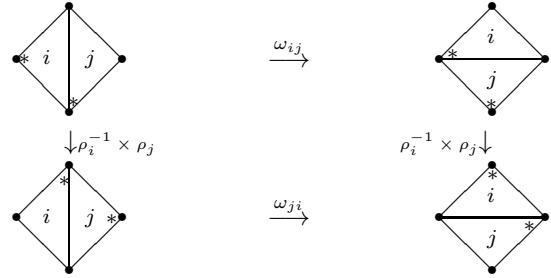
\begin{figure}[htb]
  \centering
  \begin{picture}(200,100)
\put(0,60){\begin{picture}(200,40)
\put(0,0){\begin{picture}(40,40)
\put(20,0){\line(-1,1){20}}
\put(40,20){\line(-1,-1){20}}
\put(0,20){\line(1,1){20}}
\put(40,20){\line(-1,1){20}}
\put(20,0){\line(0,1){40}}
\put(20,0){\circle*{3}}
\put(0,20){\circle*{3}}
\put(20,40){\circle*{3}}
\put(40,20){\circle*{3}}
\footnotesize
\put(10,18){$i$}\put(26,18){$j$}
\put(1,18){$*$}
\put(19.5,2){$*$}
\end{picture}}
\put(160,0){\begin{picture}(40,40)
\put(20,0){\line(-1,1){20}}
\put(40,20){\line(-1,-1){20}}
\put(0,20){\line(1,1){20}}
\put(40,20){\line(-1,1){20}}
\put(0,20){\line(1,0){40}}
\put(20,0){\circle*{3}}
\put(0,20){\circle*{3}}
\put(20,40){\circle*{3}}
\put(40,20){\circle*{3}}
\footnotesize
\put(18,26){$i$}\put(18,10){$j$}
\put(3,20){$*$}
\put(17.5,1){$*$}
\end{picture}}
\put(95,17){$\stackrel{\omega_{ij}}{\longrightarrow}$}
\end{picture}}
\put(18,46){$\downarrow$\tiny$(ij)\circ\rho_i\times\rho_j$}
\put(170,46){{\tiny$\rho_i$}$\downarrow$}
\put(0,0){\begin{picture}(200,40)
\put(0,0){\begin{picture}(40,40)
\put(20,0){\line(-1,1){20}}
\put(40,20){\line(-1,-1){20}}
\put(0,20){\line(1,1){20}}
\put(40,20){\line(-1,1){20}}
\put(20,0){\line(0,1){40}}
\put(20,0){\circle*{3}}
\put(0,20){\circle*{3}}
\put(20,40){\circle*{3}}
\put(40,20){\circle*{3}}
\footnotesize
\put(10,18){$j$}\put(26,18){$i$}
\put(16,2){$*$}
\put(34,17.5){$*$}
\end{picture}}
\put(160,0){\begin{picture}(40,40)
\put(20,0){\line(-1,1){20}}
\put(40,20){\line(-1,-1){20}}
\put(0,20){\line(1,1){20}}
\put(40,20){\line(-1,1){20}}
\put(0,20){\line(1,0){40}}
\put(20,0){\circle*{3}}
\put(0,20){\circle*{3}}
\put(20,40){\circle*{3}}
\put(40,20){\circle*{3}}
\footnotesize
\put(18,26){$i$}\put(18,10){$j$}
\put(32,20){$*$}
\put(17.5,1){$*$}
\end{picture}}
\put(95,17){$\stackrel{\omega_{ji}}{\longleftarrow}$}
\end{picture}}
\end{picture}
  \caption{Proof of the inversion relation~\eqref{eq:20}.}
  \label{fig:inv-rel}
\end{figure}
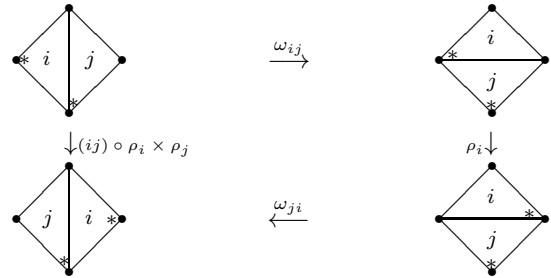

Any two (decorated ideal) triangulations can be transformed to each
other by a (finite) composition of elementary transformations, ie
$\rho_i$, $\omega_{ij}$ and permutations. This follows from the known
fact that any two ideal triangulations (without decoration) can be
transformed one into another by a composition of flips \cite{mosh},
and that
all possible decorations of a fixed ideal triangulation are
transitively acted upon by compositions of permutations and $\rho_i$ transformations.

\subsection{Quantum theory}
\label{sec:quantum-theory}

Suppose we are given a solution to eqns~\eqref{eq:1} --~\eqref{eq:4}. 
For each $\tau\in\SDIT$ and $1\le i\le \NTR$ assign
\begin{equation}
  \label{eq:17}
     \FUNCTOR(\tau,\rho_i(\tau))\equiv\ROTATE_i\in
\End\vsp^{\NTR}. 
\end{equation}
Let $i\ne j$ be such that triangles
$\bar\tau(i)$ and $\bar\tau(j)$ are as in
the lhs of figure~\ref{fe}. Then we put
    \begin{equation}\label{tij}
\FUNCTOR(\tau,\omega_{ij}(\tau))\equiv\PTOLEMY_{ij}\in
\End\vsp^{\NTR}.
    \end{equation}
Finally, for any permutation
$\sigma\in\PGROUP_{\NTR}$ set
    \[
\FUNCTOR(\tau,\tau\sigma)\equiv\PERMUTE_\sigma\in
\End \vsp^{\NTR},
    \]
    where operator $\PERMUTE_\sigma$ is defined by
    eqn~(\ref{eq:perm}).  Ensured by the consistency of
    eqns~\eqref{eq:1} --~\eqref{eq:4} with eqns~\eqref{eq:19}
    --~\eqref{eq:20}, mapping $\FUNCTOR$ can be extended to an
    operator valued function $\FUNCTOR(\tau,\tau')$ on
    $\SDIT\times\SDIT$ such that for any $\tau,\tau',\tau''\in\SDIT$
\begin{equation}
  \label{eq:10}
 \FUNCTOR(\tau,\tau)=1,\quad \FUNCTOR(\tau,\tau')\FUNCTOR(\tau',\tau'')\FUNCTOR(\tau'',\tau)\in
\COMPLEXS-\{0\}. 
\end{equation}
In particular (when $\tau''=\tau$),
\begin{equation}
  \label{eq:18}
  \FUNCTOR(\tau,\tau')\FUNCTOR(\tau',\tau)\in\COMPLEXS-\{0\}.
\end{equation}
As an example, deduce operator
$\FUNCTOR(\tau,\omega_{ij}^{-1}(\tau))$. Denoting $\tau'\equiv
\omega_{ij}^{-1}(\tau)$ and employing eqn~\eqref{eq:18} as well as
definition~\eqref{tij}, we have
\begin{equation}
  \label{eq:11}
 \FUNCTOR(\tau,\omega_{ij}^{-1}(\tau))=
\FUNCTOR(\omega_{ij}(\tau'),\tau')\simeq
(\FUNCTOR(\tau',\omega_{ij}(\tau')))^{-1}=\PTOLEMY_{ij}^{-1},
\end{equation}
where we denote by $\simeq$ an equality up to a numerical factor.

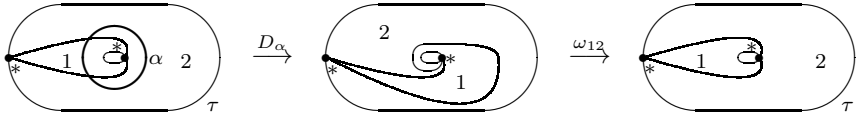
\begin{figure}[htb] 
\centering
\begin{picture}(320,40)
\put(0,0){
\begin{picture}(80,40)
\put(40,20){\oval(8,4)}
\put(40,20){\oval(80,40)}
\put(0,20){\circle*{3}}
\put(44,20){\circle*{3}}
\qbezier(0,20)(50,5)(44,20)
\qbezier(0,20)(50,35)(44,20)

\footnotesize
\put(0.5,14){$*$}
\put(39,22){$*$}
\put(20,16){$1$}
\put(65,16){$2$}
{\thicklines
\put(40,20){\circle{25}}
}
\put(53,17.5){$\alpha$}

\put(75,0){$\tau$}
\end{picture}
}
\put(95,18){$\stackrel{D_\alpha}{\longrightarrow}$}
\put(120,0){
\begin{picture}(80,40)
\put(40,20){\oval(8,4)}
\put(40,20){\oval(80,40)}
\put(0,20){\circle*{3}}
\put(44,20){\circle*{3}}
\put(38.5,20){\oval(11,10)[b]}
\put(38.5,20){\oval(11,10)[lt]}
\qbezier(38.5,25)(65,26)(65,20)
\qbezier(0,20)(68,-15)(65,20)
\qbezier(0,20)(50,5)(44,20)

\footnotesize
\put(45,17.5){$*$}
\put(0.5,13.5){$*$}
\put(20,27){$2$}
\put(49,8){$1$}

\end{picture}
}
\put(215,18){$\stackrel{\omega_{12}}{\longrightarrow}$}
\put(240,0){
\begin{picture}(80,40)
\put(40,20){\oval(8,4)}
\put(40,20){\oval(80,40)}
\put(0,20){\circle*{3}}
\put(44,20){\circle*{3}}
\qbezier(0,20)(50,5)(44,20)
\qbezier(0,20)(50,35)(44,20)

\footnotesize
\put(0.5,14){$*$}
\put(39,22){$*$}
\put(20,16){$1$}
\put(65,16){$2$}

\put(75,0){$\tau$}
\end{picture}
}
\end{picture}
\caption{The Dehn twist along contour $\alpha$ followed by
  the flip transformation $\omega_{12}$ does not change the initial
  triangulation.}\label{fannulus}
\end{figure}

The mapping class or modular group\footnote{The group of
  homeomorphisms identical on the boundary and permuting the interior
  marked points, factored wrt the connected component of the identical
  homeomorphism.}  $\MCG$ of
$\SURFACE$ naturally acts in $\SDIT$. By construction we have the
following invariance property of function $\FUNCTOR$:
\[
\FUNCTOR(f(\tau),f(\tau'))= \FUNCTOR(\tau,\tau'),\qquad \forall
f\in\MCG.
\]
This enables us to construct a projective representation of $\MCG$:
    \[
    \MCG\ni f\mapsto\FUNCTOR(\tau,f(\tau))\in
\End \vsp^{\NTR}.
    \]
    Indeed,
\[
\FUNCTOR(\tau,f(\tau))\FUNCTOR(\tau,h(\tau))=
\FUNCTOR(\tau,f(\tau))\FUNCTOR(f(\tau),f(h(\tau)))
\simeq\FUNCTOR(\tau,fh(\tau)).
\]
Any
Dehn twist, at least if it is along a non-separating contour, is
equivalent to the Dehn twist of an annulus (with two marked points on
the two boundary components) along its only
non-contractible loop denoted $\alpha$ in figure~\ref{fannulus}. As
an operator it is given in terms of operator $\PTOLEMY$.
Indeed, from figure~\ref{fannulus} it follows that $\omega_{12}\circ
D_\alpha(\tau)=\tau$.  So, using eqn~\eqref{eq:11}, we obtain
\begin{equation}\label{t01}
\FUNCTOR(\tau,D_\alpha(\tau))=\FUNCTOR(\tau,\omega_{12}^{-1}(\tau))
\simeq\PTOLEMY_{12}^{-1},
\end{equation}
where the normalization is to be fixed.

\section{Diagonalizing the Dehn twist}
\label{sec:diag-dehn-twist}
In this section we work within the quantum Teichm\"uller
theory, ie we consider
solution~\eqref{eq:5},~\eqref{eq:t-in-t-of-psi} of
system~\eqref{eq:1}~--~\eqref{eq:4}. We assume that
$(1-|\la|)\Im\la=0$,
so the range of mapping $\FUNCTOR$ in this case is given
by unitary operators.

Properly normalized quantum Dehn twist of the annulus has the form:
\[
\DEHN_{\alpha}=\zeta^{-6}\PTOLEMY_{12}^{-1}
e^{\IMUN2\pi\CONSTRAINT_\alpha^2},\quad
\CONSTRAINT_\alpha=(\MOM_1+\POS_2)/2,\quad \zeta=e^{\IMUN\pi\cla^2/3},
\]
where the exponential factor removes the spurious degree of freedom,
see~\cite{kash1,kash2}. 
Using formula~\eqref{eq:t-in-t-of-psi} we rewrite it equivalently
\[
\DEHN_{\alpha}=e^{\IMUN2\pi(\POS_\alpha^2-\cla^2)}
\QDILOG(\MOM_\alpha+\POS_\alpha),
\]
where the self-adjoint operators
\[
\MOM_\alpha\equiv\POS_1+\MOM_2-(\POS_2+\MOM_1)/2,\quad\POS_\alpha\equiv(\POS_2-\MOM_1)/2
\]
satisfy the Heisenberg commutation relation
\[
[\MOM_\alpha,\POS_\alpha]=(2\pi\IMUN)^{-1}.
\]
To diagonalize $\DEHN_\alpha$, we use the fact that mutually commuting
(and conjugate to each other in the case when $|\la|=1$) unbounded 
operators
\[
\LENGTH^\pm_\alpha\equiv
2\cosh(2\pi\la^{\pm1}\POS_\alpha)+e^{2\pi\la^{\pm1}\MOM_\alpha}
\]
 also commute with $\DEHN_\alpha$:
\[
[\LENGTH^-_\alpha,\LENGTH^+_\alpha]=[\DEHN_\alpha,\LENGTH^\pm_\alpha]=0.
\]
Geometrically operator $\LENGTH^+_\alpha$ is nothing else but the
quantized geodesic length (hyperbolic cosine thereof) of contour
$\alpha$, see~\cite{fock1,fock2}.

Let us consider the ``coordinate'' basis $\langle x|$, $x\in\REALS$,
 where $\POS_\alpha$ is diagonal
and $\MOM_\alpha$ is a differentiation operator,
\[
\langle x| \POS_\alpha =x\langle x|,\quad\langle x| \MOM_\alpha
=\frac1{2\pi\IMUN}\frac\partial{\partial x}\langle x|,
\]
and define one parameter family of vectors:
\begin{equation}
  \label{eq:6}
 \langle x|\alpha_s\rangle=\frac{\QDILOG(s+x+\cla-\IMUN
  0)}{\QDILOG(s-x-\cla+\IMUN 0)}e^{-\IMUN 2\pi(x+\cla)s}=\langle
x|\alpha_{-s}\rangle,\quad s\in\REALS .
\end{equation}

\begin{theorem}\label{theor:spec}
  Vectors $|\alpha_s\rangle$ are eigenvectors of
  operators $\LENGTH^\pm_\alpha,\DEHN_\alpha$:
  \begin{equation}
    \label{eq:7}
   \LENGTH^\pm_\alpha|\alpha_s\rangle=|\alpha_s\rangle
2\cosh(2\pi\la^{\pm1} s),\quad
\DEHN_\alpha|\alpha_s\rangle=|\alpha_s\rangle e^{\IMUN 2\pi
  (s^2-\cla^2)}. 
  \end{equation}
They are orthogonal to each other:
\begin{equation}
  \label{eq:8}
\langle\alpha_r|\alpha_s\rangle=\nu(s)^{-1} \delta(r-s),\quad
\nu(s)\equiv 4\sinh(2\pi\la s)\sinh(2\pi\la^{-1} s),  
\end{equation}
and complete in $L^2(\REALS)$:
\begin{equation}
  \label{eq:9}
\int_0^\infty|\alpha_s\rangle\nu(s)\,ds\langle \alpha_s|=1.  
\end{equation}
\end{theorem}
The continuous spectrum of Dehn twists in quantum Teichm\"uller theory
has also been  observed by V.~Fock\footnote{Private communication.}.
It is worth noting that measure $\nu(s)ds$ in
eqn~\eqref{eq:9} also appears in the representation theory of the
non-compact quantum group $\mathcal
U_q(\mathfrak{sl}(2,\REALS))$ \cite{pon1,pon2}. In the next section we
reveal that this is not accidental.

\section{Braiding and $R$-matrix}
\label{sec:braiding-r-matrix}
Braiding of a disk with two vertices is naturally associated with the 
contour (isotopy class thereof) surrounding the vertices. Moreover,
square of the braiding is Dehn twist along the
associated contour.

\begin{figure}[ht] \centering
\begin{picture}(320,180)
\put(0,100){
\begin{picture}(320,0)

\put(0,0){
\begin{picture}(80,80)
\put(40,40){\oval(80,80)}
\put(60,40){\line(-1,2){20}}
\put(40,0){\line(-1,2){20}}

\put(20,40){\line(1,2){20}}
\put(40,0){\line(1,2){20}}
\put(40,0){\line(0,1){80}}
\put(40,0){\circle*{3}}
\put(40,80){\circle*{3}}
\put(20,40){\circle*{3}}
\put(60,40){\circle*{3}}

\footnotesize
{\thicklines
\put(40,40){\oval(55,15)}}
\put(64,47){$\alpha$}

\put(14,37){$*$}\put(21,37){$*$}
\put(54,37){$*$}\put(61,37){$*$}

\put(4,36){$1$}\put(30,36){$2$}
\put(44,36){$3$}\put(71,36){$4$}

\put(75,0){$\tau$}
\end{picture}}

\put(95,38){$\stackrel{B_\alpha^{-1}}{\longrightarrow}$}

\put(120,0){
\begin{picture}(80,80)
\put(40,40){\oval(80,80)}
\put(60,40){\line(-1,2){20}}
\put(40,0){\line(-1,2){20}}
\put(40,0){\circle*{3}}
\put(40,80){\circle*{3}}
\put(20,40){\circle*{3}}
\put(60,40){\circle*{3}}
\qbezier(40,0)(-17,67)(40,40)
\qbezier(40,80)(97,13)(40,40)
\qbezier(40,0)(-45,95)(60,40)
\qbezier(40,80)(125,-15)(20,40)

\footnotesize
\put(16,40){$*$}\put(22,32){$*$}
\put(52,42){$*$}\put(60,35){$*$}

\put(30,63){$1$}\put(8,50){$2$}
\put(68,23){$3$}\put(47,10){$4$}
\end{picture}}

\put(215,38){$\stackrel{\rho_3^{-1}}{\longrightarrow}$}

\put(240,0){
\begin{picture}(80,80)
\put(40,40){\oval(80,80)}
\put(60,40){\line(-1,2){20}}
\put(40,0){\line(-1,2){20}}
\put(40,0){\circle*{3}}
\put(40,80){\circle*{3}}
\put(20,40){\circle*{3}}
\put(60,40){\circle*{3}}
\qbezier(40,0)(-17,67)(40,40)
\qbezier(40,80)(97,13)(40,40)
\qbezier(40,0)(-45,95)(60,40)
\qbezier(40,80)(125,-15)(20,40)
\footnotesize
\put(25,16){$*$}\put(22,32){$*$}
\put(52,42){$*$}\put(60,35){$*$}

\put(30,63){$1$}\put(8,50){$2$}
\put(68,23){$3$}\put(47,10){$4$}
\end{picture}}
\end{picture}}

\put(285,86){$\downarrow$ \scriptsize$\rho_1\times\omega_{23}$}
\put(45,86){$\uparrow$\scriptsize$(13)(24)\circ(\rho_3\times\rho_1^{-1})$}
\put(0,0){
\begin{picture}(320,0)

\put(0,0){
\begin{picture}(80,80)
\put(40,40){\oval(80,80)}
\put(60,40){\line(-1,2){20}}
\put(40,0){\line(-1,2){20}}

\put(20,40){\line(1,2){20}}
\put(40,0){\line(1,2){20}}
\put(40,0){\line(0,1){80}}
\put(40,0){\circle*{3}}
\put(40,80){\circle*{3}}
\put(20,40){\circle*{3}}
\put(60,40){\circle*{3}}

\footnotesize
\put(32,0){$*$}\put(21,37){$*$}
\put(40,66){$*$}\put(61,37){$*$}

\put(8,36){$3$}\put(30,36){$4$}
\put(46,36){$1$}\put(70,36){$2$}
\end{picture}}

\put(95,38){$\stackrel{\omega_{14}}{\longleftarrow}$}

\put(120,0){
\begin{picture}(80,80)
\put(40,40){\oval(80,80)}
\put(60,40){\line(-1,2){20}}
\put(40,0){\line(-1,2){20}}
\put(40,0){\circle*{3}}
\put(40,80){\circle*{3}}
\put(20,40){\circle*{3}}
\put(60,40){\circle*{3}}
\put(20,40){\line(1,0){40}}
\put(20,40){\line(1,2){20}}
\put(40,0){\line(1,2){20}}

\footnotesize
\put(32,0){$*$}\put(22,35){$*$}
\put(38,71){$*$}\put(60,35){$*$}

\put(37,53){$1$}\put(10,53){$3$}
\put(65,20){$2$}\put(37,20){$4$}
\end{picture}}

\put(210,38){$\stackrel{\omega_{13}\times\omega_{24}}{\longleftarrow}$}

\put(240,0){
\begin{picture}(80,80)
\put(40,40){\oval(80,80)}
\put(60,40){\line(-1,2){20}}
\put(40,0){\line(-1,2){20}}
\put(40,0){\circle*{3}}
\put(40,80){\circle*{3}}
\put(20,40){\circle*{3}}
\put(60,40){\circle*{3}}
\put(20,40){\line(1,0){40}}
\qbezier(40,0)(-45,95)(60,40)
\qbezier(40,80)(125,-15)(20,40)

\footnotesize
\put(25,16){$*$}\put(22,32){$*$}
\put(36,72){$*$}\put(60,35){$*$}

\put(30,63){$1$}\put(10,46){$3$}
\put(66,25){$2$}\put(47,10){$4$}
\end{picture}}
\end{picture}}

 \end{picture}
\caption{Braiding along contour $\alpha$ followed by the indicated
  sequence of transformations brings back to the initial
  triangulation $\tau$.}\label{braiding}
\end{figure}
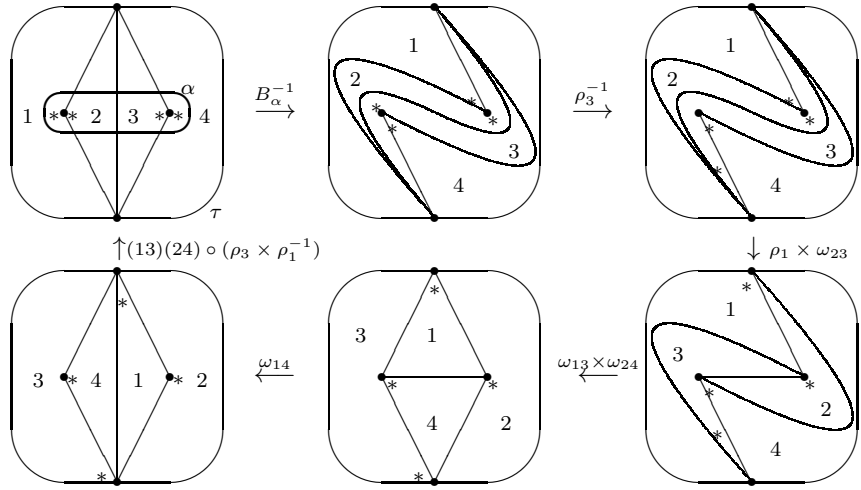

Figure~\ref{braiding} shows how the braiding of triangulation $\tau$ of a
disk with
two interior and two boundary marked points
 can be removed by a sequence of
elementary transformations:
\[
\tau=(13)(24)\circ(\rho_3\times\rho_1^{-1})\circ\omega_{14}\circ
(\omega_{13}\times\omega_{24})\circ(\rho_1\times\omega_{23})\circ
\rho_3^{-1}\circ B_\alpha^{-1}(\tau).
\]
Using the construction of subsection~\ref{sec:quantum-theory}, the
corresponding quantum braiding operator has the form 
(up to a normalization factor)
\[
\BRAID_\alpha\equiv
\FUNCTOR(\tau,B_\alpha(\tau))\simeq \PERMUTE_{(13)(24)}\RMAT_{1234},
\]
with
\begin{equation}
  \label{eq:12}
 \RMAT=\RMAT_{1234}\equiv\ROTATE_1^{-1}
\ROTATE_3\PTOLEMY_{41}\PTOLEMY_{31}\PTOLEMY_{42}
\PTOLEMY_{32}\ROTATE_1\ROTATE_3^{-1} =\PTOLEMY_{1\check4}
\PTOLEMY_{13}\PTOLEMY_{42}
\PTOLEMY_{\hat32},
\end{equation}
where in the second equation we have used notation~\eqref{eq:14} and
eqn~\eqref{eq:15}.
As a consequence of relations~\eqref{eq:1} --~\eqref{eq:4}
operator $\RMAT\in \End\vsp^4$ solves the following
Yang--Baxter equation
\begin{equation}
  \label{eq:13}
  \RMAT_{1234}\RMAT_{1256}\RMAT_{3456}=\RMAT_{3456}\RMAT_{1256}
\RMAT_{1234},
\end{equation}
and in the case corresponding to realization through the quantum
dilogarithm it is in fact (as will be shown below) $\mathcal
U_q(\mathfrak{sl}(2))$ (more precisely, the modular double thereof
\cite{fad99}) universal $R$-matrix evaluated at a reducible infinite
dimensional representation acting in $L^2(\REALS^2)$. Note that
eqn~\eqref{eq:12} is a particular case of the realization of universal
$R$-matrices in the Drinfeld doubles of Hopf algebras in terms of
solutions of the pentagon equation associated with the Heisenberg
doubles \cite{kash0}.

\subsection{Relation to Hopf algebras }
\label{sec:relat-hopf-algebr}

To identify the generators of the quantum group, we first rewrite the
$R$-matrix in the form
\begin{equation}
  \label{eq:21}
 \RMAT=\Ad(\PTOLEMY_{1\hat2}^{-1}\PTOLEMY_{\check43})\PTOLEMY_{\hat2\check4}=
\Ad(\ROTATE_2\PTOLEMY_{12}^{-1}\ROTATE_4^{-1}\PTOLEMY_{43})\PTOLEMY_{24}, 
\end{equation}
where
\[
\Ad(\a)\x\equiv \a\x\a^{-1},
\]
and we have used the pentagon equation~\eqref{eq:2}.  Let now
operator $\PTOLEMY$ be written as a sum
\[
\PTOLEMY=\sum_a \e_a\otimes\e^a,
\]
where two operator sets $\{\e_a\}$ and $\{\e^a\}$ are realizations of
mutually dual linear bases of two (mutually dual) Hopf sub-algebras
$A_\PTOLEMY$ and $A_\PTOLEMY^*$ in the associated to $\PTOLEMY$
Heisenberg double $H(A_\PTOLEMY)$.  Their co-products are given by the
formulae
\[
\Delta(\e_a)=\Ad(\PTOLEMY^{-1})(1\otimes\e_a),\quad
\Delta(\e^a)=\Ad(\PTOLEMY)(\e^a\otimes1).
\]
From eqn~\eqref{eq:21} we conclude that the $R$-matrix is also
decomposed into a sum
\[
\RMAT=\sum_a \E_a\otimes\E^a
\]
with
\begin{gather}
  \label{eq:24}
  \E_a=\Ad(\ROTATE_2\PTOLEMY^{-1})(1\otimes\e_a)=
  \Ad(\ROTATE_2)\Delta(\e_a),\\\label{eq:25}
  \E^a=\Ad(\ROTATE_2^{-1}\PTOLEMY_{21})(1\otimes\e^a)=\Ad(\ROTATE_2^{-1})\Delta'(\e^a),
\end{gather}
where $\Delta'$ is the opposite co-product. These are also
realizations of the Hopf algebras $A_\PTOLEMY$ and $A_\PTOLEMY^*$, but
this time within the Drinfeld double $D(A_\PTOLEMY)$.

\subsection{Quantum Teichm\"uller theory and $\mathcal
  U_q(\mathfrak{sl}(2))$}
\label{sec:quant-teichm}

Here we will see what quantum group and in which representation does
correspond to our solution~\eqref{eq:5}, \eqref{eq:t-in-t-of-psi}.

\begin{figure}[htb]
  \centering
  \begin{picture}(80,40)
\put(40,20){\oval(80,40)}
\put(40,20){\circle*{3}}
\put(40,40){\circle*{3}}
\put(40,0){\circle*{3}}
\put(40,0){\line(0,1){40}}
\thicklines
\put(40,20){\circle{20}}
\footnotesize
\put(34,17){$*$}
\put(41.5,17){$*$}
\put(13,17){$1$}
\put(61,17){$2$}
\put(47,27){$\beta$}
\end{picture}
  \caption{Triangulated disk with one interior vertex.}
  \label{fig:ver}
\end{figure}
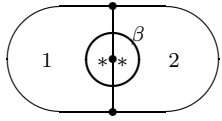

From eqns~\eqref{eq:t-in-t-of-psi},~\eqref{eq:ratio} we see that all
elements $\{\e_a\}$ can be thought to be generated by operators $\MOM$
and $e^{2\pi\la^{\pm1}\POS}$ with the co-product
\[
\Delta(\MOM)=\MOM_1+\MOM_2,\quad
\Delta(e^{2\pi\la^{\pm1}\POS})=e^{2\pi\la^{\pm1}(\POS_1+\MOM_2)}+e^{2\pi\la^{\pm1}\POS_2}.
\]
Similarly dual elements $\{\e^a\}$ are generated by $\POS$ and $
e^{2\pi\la^{\pm1}(\MOM-\POS)}$ with the co-product
\[
\Delta(\POS)=\POS_1+\POS_2,\quad
\Delta(e^{2\pi\la^{\pm1}(\MOM-\POS)})=
e^{2\pi\la^{\pm1}(\MOM_1-\POS_1-\POS_2)}+e^{2\pi\la^{\pm1}(\MOM_2-\POS_2)}.
\]
From eqn~\eqref{eq:24} we deduce that operators
\begin{equation}
  \label{eq:26}
\g_{12}\equiv\MOM_1-\POS_2,\quad \f_{12}^\pm\equiv
e^{2\pi\la^{\pm1}(\POS_1-\POS_2)}+ e^{2\pi\la^{\pm1}(\MOM_2-\POS_2)}  
\end{equation}
generate the set $\{\E_a\}$, while according to eqn~\eqref{eq:25} the
generators of the dual set $\{\E^a\}$ are the operators
\begin{equation}
  \label{eq:27}
\POS_1-\MOM_2=-\g_{21},\quad
e^{2\pi\la^{\pm1}(\POS_2-\POS_1)}+e^{2\pi\la^{\pm1}(\MOM_1-\POS_1)}=
\f_{21}^\pm. 
\end{equation}
In what follows we restrict our attention to the sub-algebra
corresponding to the positive exponent of parameter $\la$.
\begin{proposition}
  Operators $\f_{mn}\equiv\f_{mn}^+$, $\g_{mn}$ ($mn=12$ or $21$) with
  the relations
\begin{gather*}
  [\g_{12},\g_{21}]=0,\quad [\g_{mn},\f_{mn}]=-\IMUN\la\f_{mn},
  \quad[\g_{nm},\f_{mn}]=\IMUN\la\f_{mn},\\
  [\f_{12},\f_{21}]=(q-q^{-1})(e^{2\pi\la\g_{12}}-e^{2\pi\la\g_{21}}),
  \quad q\equiv e^{\IMUN\pi\la^2};
\end{gather*}
and the twisted co-product
\[
\Delta_\varphi=\Ad(e^{\IMUN
  \varphi(\g_{21}\otimes\g_{12}-\g_{12}\otimes\g_{21})})\circ\Delta,
\quad\varphi\in\REALS,
\]
\begin{gather*}
  \Delta(\g_{mn})=\g_{mn}\otimes1+1\otimes\g_{mn},\\
  \Delta(\f_{12}) =\f_{12}\otimes
  e^{2\pi\la\g_{12}}+1\otimes\f_{12},\quad \Delta(\f_{21})
  =e^{2\pi\la\g_{21}}\otimes \f_{21}+\f_{21}\otimes1.
\end{gather*}
generate a Hopf algebra $\mathcal G_\varphi$.
\end{proposition}
Algebra $\mathcal G_{\pi/2}$ is closely related with $\mathcal
U_q(\mathfrak{sl}(2))$.
\begin{definition}\footnote{This definition is that of \cite{pon2}
    without concretizing the star-structure, since the latter
    depends on a solution to the equation $(1-|\la|)\Im\la=0$.}
  Quantum group $\mathcal U_q(\mathfrak{sl}(2))$ is a Hopf algebra
  with
\begin{itemize}
\item generators $E,F,K,K^{-1}$;
\item relations
\[
KE=qEK,\quad KF=q^{-1}FK,\quad [E,F]=-(K^2-K^{-2})/(q-q^{-1});
\]
\item co-product
\[
\Delta(K)=K\otimes K,\quad\Delta(X)=X\otimes K+K^{-1}\otimes X,\quad
X=E\ \mathrm{or}\ F.
\]
 \end{itemize}
\end{definition}
\begin{proposition}
  There exists a faithful Hopf algebra homomorphism
\[
\eta\colon\mathcal U_q(\mathfrak{sl}(2))\hookrightarrow\mathcal
G_{\pi/2}
\]
such that
\begin{gather*}
  K\mapsto e^{\pi\la(\g_{12}-\g_{21})/2},\\ E\mapsto
  e^{-\pi\la(\cla+\g_{21})}\frac{\f_{21}}{q-q^{-1}},\quad
  F\mapsto
 \frac{\f_{12}}{q-q^{-1}} e^{\pi\la(\cla-\g_{12})}.
\end{gather*}
\end{proposition}
The generators of algebra $\mathcal G_\varphi$ depend on only three
combinations of the Heisenberg operators $\MOM_i,\POS_i$ ($i=1,2$)
\begin{gather*}
  \CONSTRAINT_\beta\equiv(\MOM_1-\POS_1+\MOM_2-\POS_2)/2,\\
  \MOM_\beta\equiv(\MOM_1+\POS_1-\MOM_2-\POS_2)/2,\quad
  \POS_\beta\equiv(-\MOM_1+\POS_1+\MOM_2-\POS_2)/2
\end{gather*}
which satisfy the commutation relations
\[
[\CONSTRAINT_\beta,\MOM_\beta]=[\CONSTRAINT_\beta,\POS_\beta]=0,\quad
[\MOM_\beta,\POS_\beta]=(2\pi\IMUN)^{-1}.
\]
Subscript $\beta$ here refers to the associated to our operators
contour around the interior vertex of a triangulated disk (with one
interior vertex) in figure~\ref{fig:ver}. In particular, in quantum
Teichm\"uller theory operator $\CONSTRAINT_\beta$ describes the
spurious degree of freedom associated with contour $\beta$.
Substituting these definitions into eqns~\eqref{eq:26},~\eqref{eq:27}
we obtain
\begin{gather*}
  \g_{12}=\CONSTRAINT_\beta+\MOM_\beta,
  \quad\g_{21}=\CONSTRAINT_\beta-\MOM_\beta,\\
  \f_{12}=2e^{2\pi\la\POS_\beta}
  e^{\pi\la(\CONSTRAINT_\beta+\MOM_\beta-\cla)}
  \sinh(\pi\la(\CONSTRAINT_\beta-\MOM_\beta+\cla)),\\
  \f_{21}=-2e^{\pi\la(\CONSTRAINT_\beta-\MOM_\beta+\cla)}
  e^{-2\pi\la\POS_\beta}
  \sinh(\pi\la(\CONSTRAINT_\beta+\MOM_\beta+\cla)),
\end{gather*}
while $R$-matrix~\eqref{eq:12} takes the form
\begin{equation}
  \label{eq:28}
\RMAT=e^{\IMUN2\pi(\MOM_{\beta_1}+\CONSTRAINT_{\beta_1})(\MOM_{\beta_2}-\CONSTRAINT_{\beta_2})}
\Ad\left(\frac{\QDILOG(\MOM_{\beta_1}-\CONSTRAINT_{\beta_1})}
{\QDILOG(\MOM_{\beta_2}+\CONSTRAINT_{\beta_2})}\right)
(\QDILOG(\POS_{\beta_1}+\CONSTRAINT_{\beta_1}-\POS_{\beta_2}-\MOM_{\beta_2}))^{-1}.  
\end{equation}
On the eigenvectors $\langle\xi,x|$ of operators $\CONSTRAINT_\beta$
and $\POS_\beta$,
\[
\langle\xi,x|\CONSTRAINT_\beta=\xi\langle\xi,x|,\quad\langle\xi,x|\POS_\beta=x\langle\xi,x|,
\quad\langle\xi,x|\MOM_\beta=\frac1{2\pi\IMUN}\frac\partial{\partial
  x}\langle\xi,x|,
\]
algebra $\mathcal G_\varphi$ is irreducibly represented for each fixed
$\xi\in\REALS$, and these are exactly the representations from
Schm\"udgen's classification list \cite{schm} considered by Ponsot and
Teschner in \cite{pon1,pon2}.
\begin{theorem}\label{theor:repr} If $\Im\la=0$, the following equation holds
\[
\langle\xi,x|\eta(a)=\pi_{\IMUN(\xi-\cla)}(a)\langle\xi,x|,\quad
\forall a\in \mathcal U_q(\mathfrak{sl}(2)),
\]
where representations $\pi_\alpha$, $\alpha\in \IMUN(\REALS-\cla)$,
are defined in eqn~(12) of \cite{pon2}.
\end{theorem}
Now it becomes clear the relevance of the result of Ponsot and
Teschner on the decomposition of the tensor products of the
representations under consideration to the spectral problem of Dehn
twists in quantum Teichm\"uller theory considered in
section~\ref{sec:diag-dehn-twist}: the decomposition of tensor
products of representations into irreducibles is essentially
equivalent to the spectral problem of the $R$-matrix, while square of
the latter is nothing else but a Dehn twist. Elaboration of this
connection is to be published elsewhere.

\section{Conclusion}
\label{sec:conclusion}
In this paper we have described the spectrum of Dehn twists
(Theorem~\ref{theor:spec}) in the quantum Teichm\"uller theory and
demonstrated appearance of certain infinite dimensional
representations of the quantum group $\mathcal U_q(\mathfrak{sl}(2))$
(Theorem~\ref{theor:repr}) studied in \cite{schm,pon1,pon2}. This
indicates a relationship between the two approaches to quantum
Liouville theory: one through quantization of the Teichm\"uller spaces
and another through representation theory of a non-compact quantum
group.  The result of \cite{buf} on the representation theory of the
quantum Lorentz group is also likely to be relevant here. Besides,
there should exist direct connections to quantum Liouville theory on a
space-time lattice \cite{fkv} where the non-compact quantum
dilogarithm plays the central role too.
 
We have derived the $R$-matrix (associated with braidings in the
mapping class groups) in terms of the non-compact quantum dilogarithm,
formula~\eqref{eq:28}, which first has been suggested by Faddeev in
\cite{fad99} as the universal $\mathcal U_q(\mathfrak{sl}(2))$
$R$-matrix for the corresponding modular double.  Note that more
general formula~\eqref{eq:12} directly follows from the canonical
embedding of the Drinfeld doubles of Hopf algebras into tensor product
of two Heisenberg doubles \cite{kash0}.  Thus the explanation of
section~\ref{sec:braiding-r-matrix} can be considered as a geometrical
view on the pure algebraic construction of \cite{kash0}.  It is also
worth mentioning that this $R$-matrix is the direct non-compact
analogue of the finite dimensional $R$-matrix introduced in
\cite{kash3} which leads to a specialization of the colored Jones link
invariants (polynomials) \cite{mur} with the particularly remarkable
asymptotic behavior \cite{kash4}. In this light it would be
interesting to study the corresponding non-compact analogue of the
Jones invariants.

\end{document}